\input amstex
\documentstyle{amsppt}
\magnification=\magstep1 \pagewidth{6.2 in} \pageheight{7.7in}
\hcorrection{-0.4in} \vcorrection{-0.4in} \abovedisplayskip=10pt
\belowdisplayskip=10pt
\parskip=4pt
\parindent=5mm
\baselineskip=2pt
\NoBlackBoxes

\topmatter

\title On a $q$-analogue of the $p$-adic generalized twisted $L$-functions and
$p$-adic $q$-integrals
\endtitle
\author  Lee-Chae Jang and Tae-Gyun Kim\endauthor
\affil{{\it Department of Mathematics and Computater Science}\\
{\it KonKuk University, Chungju 380-701, Korea }\\
{\it e-mail:leechae-jang$\@$hanmail.net}\\\\
{\it Ju-Kong APT  103-Dong 1001-Ho}\\
{\it Young-Chang Ri  544-4,  Hapcheon-Up Hapchon-Gun Kyungshang Nam-Do, postal no. 678-802, S. Korea }\\
{\it e-mail:jany69$\@$hanmail.net}}\\
\endaffil

\define\Q{\Bbb Q_p}
\define\C{\Bbb C_p}
\define\BZ{\Bbb Z}
\define\BQ{\Bbb Q}
\define\BC{\Bbb C}
\define\BN{\Bbb N}
\define\BR{\Bbb R}
\define\Z{\Bbb Z_p}

\keywords $p$-adic integrals, $p$-adic twisted $L$-functions,
$q$-Bernoulli numbers \endkeywords \subjclass 11B68, 11S80
\endsubjclass \abstract{  The purpose of this paper is to define
generalized twisted $q$-Bernoulli numbers by using $p$-adic
$q$-integrals. Furthermore, we construct a $q$-analogue of the
$p$-adic generalized twisted $L$-functions which interpolate
generalized twisted $q$-Bernoulli numbers. This is the
generalization of Kim's $h$-extension of $p$-adic $q$-$L$-function
which was constructed in [5] and is a partial answer for the open
question which was remained in [ 3 ] .}
\endabstract
\rightheadtext{ L.C. Jang and  T. G. Kim} \leftheadtext{On a
$q$-analogue of the $p$-adic twisted L-functions $\cdots$ }
\TagsOnRight
\endtopmatter

\document

\head \S 1. Introduction \endhead

 Let us denote $\BN, \BZ,\BQ, \BR, \BC$ sets of positive integer, integer,
rational, real and complex numbers respectively. Let $p$ be  prime
and $x \in \BQ$. Then $x= p^{v(x)} {m \over n}$, where $m,n,
v=v(x) \in \BZ$, $m$ and $n$ are not divisible by $p$. Let $|x|_p
= p^{-v(x)}$ and $ |0|_p =0$. Then $|x|_p$ is valuation on $\BQ$
satisfying
$$ |x+y|_p \leq \max \{ |x|_p , |y|_p \} .$$
Completion of $\BQ$ with respect to $| \cdot |_p$ is denoted by
$\Q$ and called the field of $p$-adic rational numbers. $\C$ is
the completion of algebraic closure of $\Q$ and $\Z = \{ x \in \Q
~|~ |x|_p  \leq 1 \} $ is called the ring of $p$-adic rational
integers(see [1,2,10,12,16]).

Let $l$ be a fixed integer and let $p$ be a fixed prime number.
We set
$$\aligned
&X=\varprojlim_N (\BZ/lp^N \BZ),\\
&X^*=\bigcup\Sb 0<a<lp\\ (a,p)=1\endSb (a+lp\Z ),\\
&a+lp^N \Z =\{x\in X\mid x\equiv a\pmod{lp^N}\},
\endaligned $$
where $N \in \BN$ and  $a\in \BZ$ lies in $0\leq a<lp^N$, cf.
[3,7,8,9].

When one talks of $q$-extension, $q$ is considered in many ways
such as an indeterminate, a complex number $q \in \BC$, or a
$p$-adic number $q \in \C$. If $q \in \BC$, one normally assumes
$|q|<1$. If $q \in \C$, then we assume $|q-1|_p < p^{- {1 \over
{p-1}}}$, so that $q^x = \exp (x \log q)$ for each $x \in X$. We
use the notation as $[x] = [x; q]= {{1-q^x}\over {1-q}}$ for each
$x \in X$. Hence $ \lim_{q \to 1} [x] =x$, cf.[4, 16, 18, 19, 20].
For any positive integer $N,$ we set
$$\mu_q (a+lp^N\Z)=\frac{q^a }{[lp^N ]}, \text{ cf. [5, 6, 7, 8, 9, 10, 11, 12, 13, 14] ,}$$
and this can be extended to a distribution on $X$. This
distribution yields an integral for each nonnegative integer $n$
(see [7]) :
$$ \int_{\Z}  [x]^n \,d\mu_q(x)=\int_X  [x]^n \,d \mu_q (x)= \beta_n (q) ,$$
where $\beta_n (q)$ are the n-th Carlitz's $q$-Bernoulli number,
cf. [4, 5, 6, 7, 8, 9, 10, 11, 12, 13, 14].

In the paper [17], Koblitz constructed $p$-adic $q-L-$function
which interpolates Carlitz's $q$-Bernoulli number at non-positive
integers and suggested two questions. One of these two questions
was solved by Kim (see [7]). In fact, Kim constructed $p$-adic
$q$-integral and proved that Carlitz's $q$-Bernoulli number can be
represented as an $p$-adic $q$-integral by the $q$-analogue of the
ordinary $p$-adic invariant measure. And also Kim is constructed a
$h$-extension of $p$-adic $q-L-$function which interpolates the
$h$-extension of $q$-Bernoulli numbers at non-positive integers
(see [5, 6, 7, 6, 9, 10, 11, 12, 13, 14]). In [5, 12, 13], Kim
constructed  $p$-adic $q$-$L$-functions and he studied their
properties . In [5], Kim introduced the $h$-extension of $p$-adic
$q$-$L$-functions and investigated many interesting physical
meaning. Also,  In [15, 16],  Koblitz  defined $p$-adic twisted
$L$-functions , and  he constructed $p$-adic measures and
integrations. And also Kim et al [3] constructed a $q$-analogue of
the twisted Dirichlet's  $L$-function which interpolated the
twisted Carlitz's $q$-Bernoulli numbers, and they remained an open
question in [3] as follows:

Find $q$-analogue of the $p$-adic twisted $L$-function which
interpolates $q$-Bernoulli numbers $\beta_{m,w,\chi}^{(h)}(q)$,
 by means of a method provided by Kim, cf. [5].

 In this paper, we will construct the
 "twisted" $p$-adic generalized $q-L-$functions and
generalized $q$-Bernoulli numbers to be a part of answer for the
question which was remained by Kim et al in [3] by means of the
same method provided by Kim in [5: p.98]. In section 2, we
construct generalized twisted $q$-Bernoulli polynomials by using
$p$-adic $q$-integrals by the same method of Kim, cf. [ 3, 5, 12,
13, 14, 20, 21, 22, 23]. We prove a formula between generalized
twisted $q$-Bernoulli polynomials which is regarded as a
generalization of Witt's formula for Carlitz's $q$-Bernoulli
polynomials in [5, Eq (5.9)], [13] and [7, Theorem 2]. This means
that the $q$-analogue of generalized twisted $q$-Bernoulli numbers
occur in the coefficients of some stirling type series. We also
give construction of the distribution of the $p$-adic generalized
twisted $q$-Bernoulli distribution. In section 3, we define the
$p$-adic generalized twisted $L$-function and construct a
$q$-analogue of the $p$-adic generalized twisted $L$-function
which interpolate generalized twisted $q$-Bernoulli numbers on
$X$. This result is related as a generalization of a $q$-analogue
of the $p$-adic $L$-function which interpolate Carlitz's
$q$-bernoulli numbers in [5, 11, 12, 13], of $p$-adic generalized
$L$-function which interpolates the $h$-extension of $q$-Bernoulli
numbers at non-positive integers in [5, 6, 7].

\head \S 2. Generalized twisted $q$-Bernoulli polynomials
\endhead

In this section, we give generalized twisted $q$-Bernoulli
polynomials by using $p$-adic $q$-integrals on $X$. Let $UD(X)$ be
the set of uniformly differentiable functions on $X$. For any $f
\in UD(X)$, T. Kim defined a $q$-analogue of an integral with
respect to an $p$-adic invariant measure in [5] which is called
$p$-adic $q$-integral. The $p$-adic $q$-integral was defined as
follows:
$$ \aligned
I_q (f) & = \int_X f(x) ~d \mu_q (x)  \\
 & = \lim_{N \to \infty} {1 \over {[lp^N ]}} \sum_{0\leq x <
 lp^N } f(x) q^x ,
\endaligned \tag 1$$
cf. [4,5,6,7,8]. Note that
$$ \aligned
I_1 (f) & = \lim_{q \to 1} I_q (f)= \int_X f(x) ~d \mu_1 (x)  \\
 & = \lim_{N \to \infty} {1 \over {lp^N }} \sum_{0\leq x <
 lp^N } f(x) ,
\endaligned \tag 2$$
and that
$$ I_1 (f_1 ) = I_1 (f) + f' (x), \tag 3$$
where $f_1 (x) =f(x+1)$.

 Let $T_p = \cup_{n \geq 1} C_{p^n} = \lim_{n \to \infty} \BZ / p^n \BZ$,
where $C_{p^n } =\{ \xi \in X | ~ \xi^{p^n} =1 \}$ is the cyclic
group of order $p^n$, see [9]. For $ \xi \in T_p $, we denote by $
\phi_\xi : \Z \to \C$ the locally constant function $ x \mapsto
\xi^x $. If we take $f(x) = \phi_\xi (x) e\xi^{tx}$, then we have
that
$$ \int_X e^{tx} \phi_\xi (x) d\mu_1 (x) = {t \over {we^t -1}}
, \tag 4 $$ cf. [5,8]. It is obvious from (3) that

$$  \int_X e^{tx} \chi (x) \phi_\xi (x) d\mu_1 (x) =
{{\sum_{a=1}^l \chi (a) \phi_\xi (a) e^{at} } \over {\xi^l e^{lt} -1}}
. \tag 5 $$

Now we define the analogue of Bernoulli numbers as follows:
$$ \aligned
e^{xt} {t \over {\xi e^t -1 } }  & = \sum_{n=0}^\infty B_{n,\xi} (x) {t^n \over {n!}}\\
{{\sum_{a=1}^l \chi (a) \phi_\xi (a) e^{at} } \over {\xi^l e^{lt}
-1}}&= \sum_{n=0}^\infty B_{n,\xi, \chi} {t^n \over n!}  ,
\endaligned \tag 6$$
cf. [5,8]. By (4), (5) and (6), it is not difficult to see that
$$\int_X  x^n \phi_\xi (x) \,d \mu_1 (x)= B_{n, \xi} \tag 7$$
and
$$\int_X  \chi(x)x^n \phi_\xi (x) \,d \mu_1 (x)= B_{n, \xi ,\chi}. \tag 8 $$

 From (7) and (8) we consider twisted $q$-Bernoulli numbers
 by using $p$-adic $q$-integral on $\Bbb{Z}_p$. For $\xi \in T_p$ and $h \in \Bbb{Z},$
 we define twisted $q$-Bernoulli polynomials as
 $$ \beta^{(h)}_{m,\xi}(x,q) = \int_{\Bbb{Z}_p} q^{(h-1)y} \xi^y [x+y]^m d\mu_q(y).
 \tag 9 $$
 Observe that
 $$ \lim_{q\rightarrow 1}  \beta^{(h)}_{m,\xi}(x, q) = B_{m,\xi}(x). $$
When $x=0$, we write  $\beta^{(h)}_{m,\xi}(0,q) =
\beta^{(h)}_{m,\xi}(q)$, which are called twisted $q$-Bernoulli
numbers. It follows from (9) that
$$ \beta^{(h)}_{m,\xi}(x,q) = \frac{1}{(1-q)^{m-1}} \sum_{k=0}^m {m \choose k} q^{xk}
(-1)^k \frac{k+h}{1 - q^{h+k} \xi}. \tag 10 $$ The Eq.(10) is
equivalent to
$$\beta^{(h)}_{m,\xi}(q) = -m \sum_{n=0}^\infty [n]^{m-1}q^{hn} \xi^n -
(q-1)(m+h) \sum_{n=0}^\infty [n]^{m}q^{hn} \xi^n. \tag 11 $$ From
(9), we obtain the below distribution relation for the twisted
$q$-Bernoulli polynomials as follows. In fact, the proof of Lemma
1 is similar to the proof of Lemma 2 with $\chi =1$.

\proclaim{Lemma  1 } For $n \ge 1,$ we have
$$ \beta^{(h)}_{n,\xi}(x,q) = d^{n-1} \sum_{a=0}^{d-1} \xi^a q^{ha} \beta^{(h)}_{n,\xi^d}(\frac
ad,q^d). $$
\endproclaim

For $\xi \in T_p$ and $h \in X$, we define generalized twisted
$q$-Bernoulli polynomials as
$$  \beta_{n,\xi,\chi}^{(h)} (x, q) = \int_X \chi (y) q^{(h-1)y}
\xi^y [x+y]^n d \mu_q (y). \tag 12$$ Observe that when $\chi=1$,
$$ \beta_{n,\xi,1}^{(h)} (x, q) = \int_X  q^{(h-1)y}
\xi^y [x+y]^n d \mu_q (y)= \beta_{n,\xi}^{(h)} (x, q) \tag 13$$
and
$$ \lim_{q \to 1} \beta_{n,\xi,\chi}^{(h)} (x, q) = \int_X  \chi (y)
\xi^y [x+y]^n d \mu_1 (y)= B_{n,\xi, \chi}^{(h)} (x),  \tag 14$$
where $\beta_{n,\xi}^{(h)} (x,q) $ is a twisted $q$-Bernoulli
polynomial and $B_{n,\xi, \chi}^{(h)} (x) $ is a generalized
Bernoulli polynomial.

\proclaim{Lemma  2 }
 For any  $n \geq 1,$ we have
 $$\beta_{n, \xi, \chi}^{(h)}(x, q)= [l]^{n-1}
\sum_{a=0}^{l-1} \chi(a) \xi^a q^{ha} \beta_{n, \xi^l , \chi^l
}^{(h)} ({a+x \over l} , q^l ). \tag 15$$
 \endproclaim
\demo{Proof} For each $n \in \BN$, we have
$$ \aligned
\beta_{n, \xi, \chi}^{(h)}(x, q)&= \int_X \chi (y) q^{(h-1)y}
\xi^y [x+y]^n d \mu_q (y)\\
&= \lim_{N \to \infty} \sum_{x_1 =0 }^{lp^N -1} \chi (x_1 )
\xi^{x_1} [x+ x_1 ]^n \mu_q (x_1 + lp^N \Z)\\
&=\lim_{N \to \infty} { 1 \over  {[lp^N]}}  \sum_{x_1 =0 }^{lp^N
-1} \chi (x_1 )
\xi^{x_1} [x+ x_1 ]^n q^{x_1}\\
&=[l]^{n-1} \sum_{a=0}^{l-1} \chi(a) q^{ha} \xi^a  \lim_{N \to
\infty} { 1 \over {[p^N : q^l]}} \sum_{m =0 }^{p^N -1}
(q^l)^{(h-1)m} (\xi^l)^m [{{x+a} \over l} + m : q^l ]^n (q^l)^m \\
&=[l]^{n-1} \sum_{a=0}^{l-1} \chi(a) \xi^a q^{ha} \beta_{n, \xi^l
, \chi^l }^{(h)} ({a+x \over l} , q^l ).
\endaligned $$
\enddemo

We note that when $x=0$, we have the distribution relation for the
generalized twisted $q$-Bernoulli numbers as follows: for $n \geq
1$,
 $$\beta_{n, \xi, \chi}^{(h)}(q) = \beta_{n, \xi, \chi}^{(h)}(0, q)= [l]^{n-1}
\sum_{a=0}^{l-1} \chi(a) \xi^a q^{ha} \beta_{n, \xi^l }^{(h)} ({a
\over l} , q^l ) \tag 16$$ and that when $x=0$ and $q=1$, we have
the distribution relation for the generalized twisted Bernoulli
numbers as follows: for $n\geq 1$,
 $$\beta_{n, \xi, \chi}^{(h)}= \beta_{n, \xi, \chi}^{(h)}(1)= l^{n-1}
\sum_{a=0}^{l-1} \chi(a) \xi^a  \beta_{n, \xi^l }^{(h)} ({a \over
l}  ) \tag 17$$ and that when $x=0$ and $\chi=1$, we have the
distribution relation for the twisted $q$-Bernoulli polynomials as
follows : for $n\geq 1$,
$$\beta_{n, \xi}^{(h)}(0,q) = \beta_{n, \xi, 1}^{(h)}(q)= [l]^{n-1}
\sum_{a=0}^{l-1}  \xi^a q^{ha} \beta_{n, \xi^l }^{(h)} ({a \over
l} , q^l ). \tag 18$$

Lemma 1 and Lemma 2 are important for the construction of the
$p$-adic generalized twisted $q$-Bernoulli distribution as
follows.

\proclaim{Theorem  3 } Let $q \in \C$. For any  positive integers
$N, n$ and $l$, let $\mu_{n,\xi}^{(h)}$ be defined by
$$ \mu_{n,\xi}^{(h)} (a +lp^N \Z) = [lp^N ]^{n-1} q^{ha} \xi^a
\beta_{n, \xi^{lp^N }} ({a \over {lp^N}} , q^{lp^N } ).$$ Then
$\mu_{n,\xi}^{(h)}$ extends uniquely to a distribution on $X$.
 \endproclaim

 \demo{Proof}
It is suffices to show
$$ \sum_{i=1}^{p-1} \mu_{n,\xi}^{(h)} (a + ip^N + p^{N+1} \Z)
= \mu_{n,\xi}^{(h)} (a +p^N \Z).
$$
Indeed, Lemma 1 and the definition of $\mu_{n,\xi}^{(h)}$ imply
that
$$ \aligned
 &\sum_{i=1}^{p-1} \mu_{n,\xi}^{(h)} (a + ip^N + p^{N+1} \Z)\\
 &= \sum_{x=0}^{p-1} [p^{N+1} ]^{n-1} q^{h(a+x p^N )} \xi^{a+xp^N}
 \beta_{n, \xi^{p^{N+1}}}^{(h)}({{a+xp^N} \over {p^{N+1}}}, q^{p^{N+1}})\\
&=[p]^{n-1} q^{ha} \xi^a [p^N : q^p ]^{n-1} \sum_{x=0}^{p-1}
(q^{p^N})^{xh}  (\xi^{p^N})^x \beta_{n,
(\xi^{p^N})^{p}}^{(h)}({{{a \over {p^N}}+x} \over {p}},
(q^{p^N})^p ) \\
&=[p]^{n-1} q^{ha} \xi^a \beta_{n, \xi^{p^N} }^{(h)} ({a \over
{p^N}} , q^{p^N} )\\
&=\mu_{n,\xi}^{(h)} (a +p^N \Z).
\endaligned $$

 \enddemo

\head \S 3. A $q$-analogue of the $p$-adic twisted $L$-functions
\endhead

Let $\alpha \in X^* , \alpha \not= 1, n \geq 1$. By the definition
of $\mu_{n,\xi,\chi}^{(h)}$, we easily see :
$$ \aligned \int_X \chi (x) d \mu_{n, \xi }^{(h)} (x) &= \beta_{n,
\xi,\chi}^{(h)}(q)\\
\int_{pX} \chi (x) d \mu_{n, \xi }^{(h)} (x) &= [p]^{n-1} \chi (p)
\beta_{n,\xi^p,\chi}^{(h)} (q^p )\\
\int_X \chi (x) d \mu_{n;q^{1 \over \alpha}, \xi^{1 \over \alpha}
}^{(h)} (\alpha x) &=  \chi ({1 \over \alpha })
\beta_{n,\xi^{1 \over \alpha},\chi}^{(h)} (q^{1 \over \alpha}  )\\
\int_{pX} \chi (x) d \mu_{n;q^{1 \over \alpha}, \xi^{1 \over
\alpha} }^{(h)} (\alpha x) &= [p;q^{1 \over \alpha }]^{n-1} \chi
({p \over \alpha}) \beta_{n,\xi^{p \over \alpha},\chi}^{(h)} (q^{p
\over \alpha} ).
\endaligned  \tag 19 $$
For compact open set $U \subset X$, we define
$$ \mu_{n;q, \alpha, \xi}^{(h)} (U) =  \mu_{n;q, \xi }^{(h)}
(U)-\alpha^{-1} [ \alpha^{-1} ;q]^{n-1}  \mu_{n;q^{1 \over
\alpha}, \xi^{1 \over \alpha}}^{(h)} (U). \tag 20$$ By the
definition of $\mu_{n;q, \xi}^{(h)} $ and (19), we note that
$$  \aligned  \int_{X^*} \chi (x) d \mu_{n;q,\alpha, \xi }^{(h)} (x)
&=
\beta_{n,\xi,\chi}^{(h)}(q) - [p]^{n-1} \chi (p) \beta_{n, \xi^p, \chi} (q^p )\\
& - {1 \over \alpha} [ {1 \over \alpha }]^{n-1} \chi ( {1 \over
\alpha} ) \beta_{n,\xi^{1 \over \alpha},\chi}^{(h)} (q^{1 \over \alpha} )\\
& + {1 \over \alpha} [ {p \over \alpha }]^{n-1} \chi ( {p \over
\alpha} ) \beta_{n,\xi^{p \over \alpha },\chi}^{(h)} (q^{p \over \alpha})\\
&= ( 1- \chi^p ) ( 1- {1 \over \alpha} \chi^{ 1 \over \alpha} )
\beta_{n,\xi,\chi}^{(h)} ,
\endaligned  \tag 21 $$
where the operator $\chi^y = \chi^{y, n;q, \xi}$ on $f(q, \xi)$
defined by
$$ \chi^y f(q,\xi) = [y]^{n-1} \chi(y) f(q^y , \xi^y ) ,~~ \chi^x \chi^y =
chi^{x, n;q^y ,\xi^y} \circ \chi^{y, n;q, \xi} . $$  Let  $x \in
X$. We recall that $\{ x \}_N$ denote the least nonnegative
residue (mod $lp^N$) and that if $[x]_N = x-\{ x\}_N $, then
$[x]_N \in lp^N \Z$. Now we can define in [5] as follows:
$$ \mu_{Mazur, 1, \alpha}^{(h)} (a+ lp^N \Z)   = ( {{{1 \over \alpha}
-1}\over {h+1}} + { h \over \alpha} \cdot {{[a \alpha]_N} \over
{lp^N} } ) .$$ By the same method of Kim in [5], we easily see:
$$ \aligned & \lim_{N \to \infty} \mu_{n;q,\alpha,\xi}^{(h)} (a +
lp^N \Z )\\
&= \lim_{N\to \infty} [l]^{n-1} ((h+n) q^{(h+1)a} - h q^a ) \xi^a
( {{{1 \over \alpha} -1}\over {h+1}} + { h \over \alpha} \cdot
{{[a \alpha]_N} \over {lp^N} } ) .
\endaligned  \tag 22 $$
Thus we have
$$ \mu_{n;q,\alpha,\xi}^{(h)} (x) =  [x ]^{n-1}((h+n) q^{(h+1)x} - h q^{xh} ) \xi^x
\mu_{Mazur, 1, \alpha}^{(h)} (x)   .
 \tag 23 $$

\proclaim{Theorem  4 } $\mu_{n;q,\alpha,\xi}^{(h)}$ are bounded
$\C$-valued measure on $X$ for all $ n\geq 1$ and $ \alpha\in X^*
, \alpha \not= 1$.
 \endproclaim

Now we define $<x>=<x;q>=[x;q]/w(x)$, where $w(x) $ is the
Teichm\"{u}ller character. For $|q-1|_p < p^{- {1 \over {p-1}}}$,
we note that $<x>^{p^N } \equiv 1(mod~ p^N)$. By (21) and (23), we
have the following:
$$  \aligned  &\int_{X^*} \chi_n (x) d \mu_{n;q,\alpha \xi }^{(h)}
(x)\\
&= \int_{X^*} \chi_n (x)   [x]^{n-1}((h+n) q^{(h+1)x} - h q^{xh} )
\xi^x \mu_{Mazur, 1, \alpha}^{(h)} (x)\\
&= \int_{X^*} ((h+n) q^{(h+1)x} - h q^{xh} ) <x>^{n-1} \xi^x
\chi_1 (x)  \mu_{Mazur, 1, \alpha}^{(h)} (x)
\endaligned  \tag 24 $$
where $\chi_n (x) = \chi w^{-n} (x)$. By using (24), we can
construct a $q$-analogue of $p$-adic generalized twisted
$L$-function.

\proclaim{Definition  5 } For fixed $\alpha\in X^* , \alpha\not=
1$, we define a $h$-extension of $p$-adic generalized twisted
$L$-function as follows;
$$ L_{p,q,\xi}^{(h)} (s,\chi)= {1 \over {1-s}} \int_{X^*} ((h+1-s)q^{(h+1)x} -hq^{hx} ) \xi^x <x>^{-s} \chi_1 (x) d
\mu_{Mazur,1,\alpha}^{(h)} (x), \tag 25 $$ for $s \in X$.
\endproclaim

\proclaim{Theorem  6 } For each $s \in \Z$ and $\alpha \in X^* ,
\alpha \not= 1$, we have
$$  \aligned L_{p,q,\xi}^{(h)} (s, \chi ) =& {{1-s+h}\over {1-s}}
(q-1) {\sum_{n=1}^{\infty}}^* {{q^{nh} \xi^n w^{s-1} (n) } \over
{[n]^{s-1}}} \chi(n) ( {{{1 \over \alpha} -1}\over {h+1}} + { h
\over \alpha} \cdot {{[n\alpha]_N} \over
{lp^N} } )\\
&+  {\sum_{n=1}^{ \infty}}^* q^{hn} \xi^n [n]^{-s} w^{s-1} (n)
\chi(n)( {{{1 \over \alpha} -1}\over {h+1}} + { h \over \alpha}
\cdot {{[n\alpha]_N} \over {lp^N} } ). \endaligned \tag 26$$
 \endproclaim
 where ${\sum_{n=1}^\infty}^*$ means to sum over the rational
 integers prime to $p$ in the give range.
\demo{Proof} For each $s \in \Z$ and $x \in X^*$, we have
$$ \aligned
(h+1-s) q^{(h+1)x} - hq^{hx} &= h q^{hx} (q^x -1) + (1-s)q^x
q^{hx}\\
&= (q-1) q^{hx} [x] (h+1-s) + q^{hx}(1-s).
\endaligned $$
Thus
$$ \aligned
& { 1 \over {1-s}} \int_X^* ( (h+1 -s) q^{(h+1)x} - hq^{hx} )
\xi^x <x>^{-s} \chi_1 (x) d \mu_{Mazur, 1,\alpha }^{(h)} (x) \\
&={ 1 \over {1-s}} \int_X^*[ (q-1) q^{hx} [x] (h+1-s) +
q^{hx}(1-s)]\xi^x <x>^{-s} \chi_1 (x) d \mu_{Mazur, 1,\alpha }^{(h)} (x) \\
&={{1-s+h}\over {1-s}} (q-1) {\sum_{n=1}^{\infty}}^* {{q^{nh}
\xi^n w^{s-1} (n) } \over {[n]^{s-1}}} \chi(n)( {{{1 \over \alpha}
-1}\over {h+1}} + { h \over \alpha} \cdot {{[n\alpha]_N} \over
{lp^N} } ) \\
&+  {\sum_{n=1}^{ \infty}}^* q^{hn} \xi^n [n]^{-s} w^{s-1} (n)
\chi(n)( {{{1 \over \alpha} -1}\over {h+1}} + { h \over \alpha}
\cdot {{[n\alpha]_N} \over {lp^N} } )
\endaligned $$

\enddemo

The equation (26) with $h=s-1$ implies that
$$ L_{p,q,\xi, \alpha}^{(s-1)} (s, \chi ) = {\sum_{n=1}^{ \infty}}^* q^{(s-1)n} \xi^n [n]^{-s} w^{s-1} (n)
\chi(n)( {{{1 \over \alpha} -1}\over {s}} + { {s-1} \over \alpha}
\cdot {{[n\alpha]_N} \over {lp^N} } ) . \tag 27
$$
Finally for each positive integer $m$, we can construct a
$q$-analogue of the $p$-adic twisted $L$-function which
interpolate a generalized $q$-Bernoulli number.

\proclaim{Theorem  7 } For each $m \in \BN$ and $\alpha \in X^* ,
\alpha \not= 1$, we have
$$   L_{p,q,\xi}^{(h)} (1-m, \chi ) =  {- {1 \over m}}
(1-   \chi_m^p ) (1- {1 \over \alpha } \chi_m^{1 \over \alpha} )
w^{-m} \beta_{m, \xi, \chi}^{(h)}(q) . \tag 28$$
 \endproclaim

\demo{Proof} For each $s \in \Z$, by using (21), we have
$$  \aligned
&L_{p,q,\xi, \alpha}^{(h)} (s, \chi ) \\
&={ 1 \over {1-s}} \int_{X^*} ( (h+1 -s) q^{(h+1)x} - hq^{hx} )
\xi^x <x>^{-s} \chi_1 (x) d \mu_{Mazur, 1,\alpha }^{(h)} (x) \\
&={ 1 \over {1-s}} \int_{X^*} \chi_{1-s} (x)  d\mu_{1-s;q,
\alpha,\xi}(x)\\
&={ 1 \over {1-s}}(1-   \chi_{1-s}^p ) (1- {1 \over \alpha }
\chi_{1-s}^{1 \over \alpha} ) \beta_{n, \xi, \chi}^{(h)}(q) .
\endaligned$$
Thus
$$  \aligned
&L_{p,q,\xi }^{(h)} (1-m, \chi ) \\
&={ 1 \over m}(1-   \chi_m^p ) (1- {1 \over \alpha } \chi_m^{1
\over \alpha} ) \beta_{n, \xi, \chi}^{(h)}(q) .
\endaligned$$

\enddemo

Remark. In [5], Kim constructed the $h$-extension of $p$-adic
$q$-$L$-functions. And the question to inquire the existence of
the twisted $p$-adic $q$-$L$-functions was remained in [3]. This
is still open. By means of the method provided by Kim in [5], we
constructed the twisted $p$-adic $q$-$L$-function to be a part of
an answer for the question which was remained in [3].

\enddocument